\def\hV{{\widehat{V}}}
\def\bea{\begin{eqnarray}}
\def\ena{\end{eqnarray}}

\def\proof{{\it Proof}\quad}

\newtheorem{dfn}{Definition}
\newtheorem{prop}{Proposition}
\def\C{{\bf C}}
\def\Z{{\bf Z}}
\def\H{{\cal H}}

\def\qed{\quad$\Box$}

\documentstyle{article}
\begin{document}
\begin{center}
{\Large \bf Extended Vertex Operator Algebras and Monomial Bases
\\[8mm] }
{Boris Feigin\,$^1$ and Tetsuji Miwa\,$^2$\\[8mm]
September 1998}
\end{center}
\renewcommand{\thefootnote}{\fnsymbol{footnote}}
\footnote[0]
{Dedicated to James B. McGuire on the occasion of his 65th birthday.}
\renewcommand{\thefootnote}{\arabic{footnote}}
\footnotetext[1]{L.D. Landau Institute for Theoretical Physics,
Chernogolovka 142432, Russian Federation.}
\footnotetext[2]{Research Institute for Mathematical Sciences,
Kyoto University, Kyoto 606, Japan.}
\begin{abstract}
\noindent
We present a vertex operator algebra which is an extension of the
level $k$ vertex operator algebra for the
$\widehat{sl}_2$ conformal field theory.
We construct monomial basis of its irreducible representations.
\end{abstract}

\section{Introduction}
Recall the following well-known construction of the level-$1$ representations
of the Lie algebra $\widehat{gl}_2$. Let $\hV$ be the space of functions from
$S^1$ to $V\simeq{\C}^2$.
Then the irreducible representations $\widehat{gl}_2$
are realized in the space $\Lambda^{\frac{\infty}{2}}(\hV{})$
--- the semi-infinite
exterior power of $\hV$. There are many ways to define the space
$\Lambda^{\frac{\infty}{2}}(\hV)$. One approach is as follows. Consider
the Clifford algebra generated by the space $\hV{}\oplus\hV^*$
with the natural quadratic form. The irreducible representation of the Clifford
algebra is the direct sum
\bea\label{CLI}
\bigoplus_{i\in\Z}\Lambda^{\frac{\infty}{2}+i}(\hV).
\ena
If we choose a basis in $\hV$ then the basis in (\ref{CLI}) consists
of the semi-infinite wedge products of the basis vectors of $\hV$.

An alternative construction goes as follows.
Let $\widehat{V}_n:=V\otimes z^n{\C}[z^{-1}]$. Then $\hV=V\otimes\C[z,z^{-1}]$
is equal to the inductive limit:
$$\cdots\to\widehat{V}_{-1}\to\widehat{V}_0\to\widehat{V}_1\to\cdots$$
Note that $\widehat{V}_n$ is a graded space.
Let $\omega_n$ be an element of the highest degree in $\Lambda^2(\hV_n)$.
We can consider the sequence of embeddings:
$$\mu:\ \Lambda^0(\hV_0)\to\Lambda^2(\hV_1)\to\Lambda^4(\hV_2)\to\cdots$$
Here the map $\Lambda^s(\hV_n)\to\Lambda^{s+2}(\hV_{n+1})$ is the
composition of the map
$\Lambda^s(\hV_n)\to\Lambda^s(\hV_{n+1})$ and the product
$\Lambda^s(\hV_{n+1})\wedge\omega_{n+1}\to
\Lambda^{s+2}(\hV_{n+1})$. The inductive limit of the sequence $\mu$ is
$\Lambda^{\frac{\infty}{2}}(\hV)$.

The dual sequence
$$\mu^*:\ \Lambda^0(\hV_0)^*\leftarrow\Lambda^2(\hV_1)^*\leftarrow
\Lambda^4(\hV_2)^*\leftarrow\cdots$$ has a ``functional'' description.
Namely, let us identify the space $\hV_n^*$ with the space
$U_n=z^{-n-1}\C[z]dz\otimes(\C^2)^*$. Here the
pairing is given by the residue. Then $\Lambda^s(\hV_n)^*\cong
\Lambda^s(\hV_n^*)\cong\Lambda^s(U_n)$. The space $\Lambda^s(U_n)$
is the space of functions in $z_1,\ldots,z_s$ with values in
$({\C}^2)^*\otimes\cdots\otimes({\C}^2)^*$ (times $dz_1\cdots
dz_s$) which is

(a) skew-symmetric with respect to the permutations of $\{z_j\}$ and
components in the tensor product
$({\C}^2)^*\otimes\cdots\otimes({\C}^2)^*$;

(b) of the form $(z_1\cdots z_s)^{-n-1}P(z_1,\ldots,z_s)$
where $P$ is a vector-valued polynomial.

\noindent
Roughly speaking, the projective limit of the sequence $\mu^*$ is the space of
skew-symmetric functions in infinitely many variables.

In this paper we present a generalization of these level-$1$ constructions.
First of all, we construct a vertex operator algebra $A_k$
which play the role of the Clifford algebra. It is generated by the spaces
$V\otimes{\C}[z,z^{-1}]$ and $V^*\otimes{\C}[z,z^{-1}]$ where
$V\simeq\C^{k+1}$.
The idea of the construction is as follows. Consider the operator algebra of
the conformal field theory consisting of both currents and intertwiners.
The latter generate an algebra which is an extension of the vertex operator
algebra generated by the former.
However, it is not a vertex operator algebra because the relations among
these operators are not ``local''.
In a vertex operator algebra the operators placed at
distinct points must commute (or skew-commute). In some cases it is possible
to find the combinations of vertex operators which are ``local'' and generate
a vertex operator algebra. It gives us an ``algebraic'' extension of the
vertex operator algebra of currents.

For example let us start with the $\widehat{sl}_2$ conformal field theory
of level $1$. We have the vertex operators ${\C}^2(z)$ associated with 
$2$-dimensional representation of $sl_2$.
Consider the product of this theory and the free field theory. Fermions
in the Clifford algebra are the products of ${\C}^2(z)$ and
some primary fields (i.e., vertex operators) of the free field theory.

For higher
level we apply exactly the same construction. We consider the product of
the $\widehat{sl}_2$ conformal field theory of level $k$ and the
free field theory. The new vertex operator algebra is generated by the products
of $\widehat{sl}_2$ vertex operators with values in ${\C}^{k+1}$ and
the vertex operators of the free field theory.

We will study an analogue of the space
$\Lambda^{\frac{\infty}{2}}\widehat V$
and a monomial basis there. Here we describe the
``functional'' version of the semi-infinite construction. Return to the
level-$1$ case for the moment. Let us realize $({\C}^2)^*$ in the space of
polynomials in the variable $t$ of degree $\leq1$. Then the space $U_n$
can be identified with the space of functions $\{z^{-n-1}Q(t,z)\}$
where $Q$ is a polynomial in $t,z$ of degree $\leq1$ in $t$. Similarly,
$\Lambda^sU_n$ can be identified with the space of functions of the form
$(z_1\cdots z_s)^{-n-1}Q(t_1,z_1,\ldots,t_s,z_s)$
where $Q$ is a polynomial skew-symmetric with respect to the permutations
of pairs $(t_j,z_j)$, and of degree $\leq1$ in $t_j$.

For $k>1$,
let us introduce the space $(\Lambda^sU_n)^k$. It consists of expressions
of the form
\bea
(z_1\cdots z_s)^{-k(n+1)}R(t_1,z_1,\ldots,t_s,z_s)\label{POLR}
\ena
where $R$ is a polynomial in $(t_1,z_1,\ldots,t_s,z_s)$

(a) of degree $\leq k$ in $t_j,\ j=1,\ldots,s$;

(b) symmetric with respect to the permutations of the pairs $(t_j,z_j)$
if $k$ is even, and skew-symmetric otherwise;

(c) subject to the conditions
$$
\frac{\partial^{j_1}}{\partial t_1^{j_1}}
\frac{\partial^{j_2}}{\partial z_1^{j_2}}
R|_{t_1=t_2;\ z_1=z_2}=0\ \hbox{ for $j_1+j_2<k$}.
$$
In other words, $R$ has a zero of order $k$ if $t_1=t_2$ and $z_1=z_2$.

Thus, the space $\Lambda^sU_n$ is identified with some space of
polynomials, and
$(\Lambda^sU_n)^k$ is the linear span of
$\Lambda^sU_n\times\cdots\times\Lambda^sU_n$
($k$ times) where $\times$ denotes just the product of polynomials.
We have a projective system of spaces
\bea
(\Lambda^0U_0)^k\leftarrow(\Lambda^2U_1)^k\leftarrow
(\Lambda^4U_2)^k\leftarrow\cdots,
\ena
where the map sends the element (\ref{POLR}) in
$(\Lambda^{2n}U_n)^k$ to
$$
\frac1{k!(z_1\cdots z_{2n-2})^{k(n+1)}}
\left(\frac\partial{\partial t_{2n-1}}\right)^k
R(t_1,z_1,\ldots,t_{2n},z_{2n})|_{t_{2n-1}=t_{2n}=z_{2n-1}=z_{2n}=0}.
$$
Roughly speaking, the projective limit is a space of polynomials in infinitely
many variables with some conditions on diagonals. Its dual space
is our analogue of the space
$\Lambda^{\frac{\infty}{2}}\hV$. Let us denote it by
$(\Lambda^{\frac{\infty}{2}}\hV)^k$. Note that the
space ${\C}^{k+1}\otimes{\C}[z,z^{-1}]$ can be identified with the dual
to the projective limit
\bea
\cdots\leftarrow U^k_{-1}\leftarrow U^k_0\leftarrow U^k_1\leftarrow\cdots.
\ena
Therefore, our construction gives the semi-infinite ``power'' of the space
${\C}^{k+1}\otimes{\C}[z,z^{-1}]$.

In the second half of the paper, we construct a monomial basis of
$\bigl(\Lambda^{\frac{\infty}{2}}\widehat V\bigr)^k$.

An irreducible representations of the vertex operator algebra $A_k$
is a direct sum of irreducible representations of $\widehat{gl}_2$.
The latter are of the form $\pi_j\otimes{\cal H}_p$ where
$\pi_j$ is the irreducible representation for $\widehat{sl}_2$
of level $k$ and spin ${j\over2}$,
and ${\cal H}_p$ is a bosonic Fock space. The value $p$
of the zero-mode is chosen suitably. We remark that the space
$(\Lambda^{\frac{\infty}2}\hV)^k$
discussed above is $\pi_0\otimes{\cal H}_0$.
The generators of $A_k$,
which we denote by $\varphi_a(z)$ and $\varphi^*_a(z)$, act as
follows:
\bea
\varphi_a(z):\pi_j\otimes{\cal H}_p\rightleftharpoons
\pi_{k-j}\otimes{\cal H}_{p+\sqrt{k\over2}}:\varphi^*_a(z).
\ena
To be precise, the Fourier coefficients of the vertex operators
act as above.

Let $\omega$ be the highest weight vector of one of the
subspaces $\pi_j\otimes{\cal H}_p$. The Fourier coefficients
of $\varphi_a(z)$, which we denote by $\varphi_{a,n}$,
create a set of vectors in the total representation space.
The vector $\omega$ itself is created from
another highest weight vector, say $\omega'$,
by a Fourier coefficient, and we can go further back to
$\omega''$, etc. The spaces generated from $\omega,\omega',\ldots$
by the $\varphi_{a,n}$,
are increasing, and in fact, exhaust the whole space.

We prove this in three steps. In the first step,
we write the quadratic relations satisfied by the vertex operators.
We show that the space generated from $\omega$ is spanned by
a certain set of vectors, ``normal-ordered'' monomials of
$\varphi_{a,n}$ acting on $\omega$.
In the second step we show that the set of normal-ordered
monomials is linearly independent
by showing the non-degeneracy of the dual coupling. Finally,
we show that the union of the subspaces generated from $\omega,\omega',\ldots$
is equal to the whole space by calculating the characters.

Before passing we mention briefly some references
closely related to this work.

The construction of the level-$1$ vertex operators goes back to the papers
\cite{LW},\cite{KKLW},\cite{FK},\cite{DJKM},\cite{JM}.
The idea of semi-infinite construction used in this paper is
originally developed in \cite{FS} for the current generators
of $\widehat{sl}_2$.

Our normal-ordered monomials are labeled by the ``paths''
known in the solvable lattice models.
In \cite{DJKMO} paths are used to label a basis
for the higher level representations of $\widehat{sl}_r$.
The construction in that paper uses the Chevalley generators
of $\widehat{sl}_r$ in order to create the basis vectors.
This is the point of difference from the present paper.
We construct a basis of the level $k$ irreducible $\widehat{gl}_2$
modules by using the Fourier components of the vertex operators
taking values in $\C^{k+1}$ modified with bosonic vertex operators.

In the $q$-deformed situation, similar constructions
were given in \cite{KMS},\cite{KMPY}. However, their construction
does not recover the construction in this paper in the limit $q=1$.
The difference lie in the following point. The choice of the bosonic
vertex operators in our construction is uniquely determined by the
locality condition as explained. On the other hand, the choice in
\cite{KMS},\cite{KMPY} is such that the quadratic relations
are of finite forms in terms of the Fourier coefficients. These two
conditions are not compatible.

{\bf Acknowledgement}
We thank Murray Batchelor for giving us the opportunity of
writing this paper for this volume. We also thank Masaki Kashiwara for
useful discussions.
\section{Vertex operator algebra $A_k$}
In this section we construct the vertex operator algebra $A_k$
by using the vertex operator algebra of the $\widehat{sl}_2$
conformal field theory of level $k$ and that of the free bosons.

\subsection{Definition of $A_k$}
Recall first some well-known facts about minimal conformal field theories
associated with the affine Lie algebra
\bea
\widehat{sl}_2
=sl_2\otimes\C[t,t^{-1}]\oplus\C c\oplus\C d.
\ena
We identify $sl_2\otimes1\subset\widehat{sl}_2$ with $sl_2$.
We use the basis of $sl_2$:
\bea
E=\pmatrix{&1\cr0&\cr},H=\pmatrix{1&\cr&-1\cr},F=\pmatrix{&0\cr1&\cr}.
\ena
We set $X_i=X\otimes t^i$ for $X\in sl_2$.
We also use the Chevalley generators
\bea
e_0=F_1,h_0=-H_0+c,f_0=E_{-1},e_1=E_0,h_1=H_0,f_1=F_0.
\ena

Let $P=\C\Lambda_0\oplus\C\Lambda_1\oplus\C\delta$ be the weight lattice.
The dual lattice is $P^*=\C h_0\oplus\C h_1\oplus\C d$ where
$(\Lambda_0,\Lambda_1,\delta)$ and $(h_0,h_1,d)$ are dual to each other.

We consider level $k\in\Z_{\geq0}$ representations, i.e., $c=k$.
There are $k+1$ integrable highest weight representations
$\pi_0,\pi_1,\ldots,\pi_k$ of level $k$.
The representation $\pi_j$ is generated by the highest weight vector
$|j\rangle$ satisfying
$E_i|j\rangle=0$ $(i\geq0)$, $H_i|j\rangle=F_i|j\rangle=0$ $(i\geq1)$
$H_0|j\rangle=j|j\rangle$ and $d|j\rangle=0$.
These representations constitute the Verlinde
algebra. We will need the following relations in the Verlinde algebra.
\bea
\pi_k\cdot\pi_j=\pi_{k-j}.\label{VER}
\ena

To each representation $\pi_j$ we can correspond the vertex operators.
In this paper we consider the vertex operators corresponding to $\pi_k$.
Let $V$ be the $(k+1)$-dimensional representation of $sl_2$.
Fix the standard basis $\psi_0,\psi_1,\ldots,\psi_k$ in $V$:
$E\psi_a=(k-a+1)\psi_{a-1}$,
$H\psi_a=(k-2a)\psi_a$,
$F\psi_a=(a+1)\psi_{a+1}$.
We have $(k+1)$ vertex operators $\psi_a(z)=\sum_{n\in{\bf Z}}\psi_{a,n}z^{-n}$.
These are a collection of operators, and satisfy
\bea
[X_i,\psi_a(z)]&=&z^i(X\psi_a)(z)\hbox{ for $X\in\widehat{sl}_2$},
\\[0pt]
[d,\psi_a(z)]&=&-z{d\over dz}\psi_a(z).\label{DCOM}
\ena
From (\ref{VER}) it follows that they are acting from $\pi_j$ to $\pi_{k-j}$.
We fix the normalization of $\psi_a(z)$ by the condition
\bea
\langle k-j|\psi_j(z)|j\rangle=1.\label{NOR}
\ena

We use the currents
$E(z)=\sum_iE_iz^{-i-1}$,
$H(z)=\sum_iH_iz^{-i-1}$,
$F(z)=\sum_iF_iz^{-i-1}$.
We also use symbolic notations $\int H_{>0}(z)=-\sum_{i>0}{H_i\over i}z^{-i}$
and $\int H_{<0}(z)=-\sum_{i<0}{H_i\over i}z^{-i}$.
The currents act on the integrable representations of level $k$ and satisfy
\bea\label{INT}
E(z)^{k+1}=F(z)^{k+1}=0.
\ena
The semi-infinite construction given in \cite{FS} is based on (\ref{INT}).

Consider the group element
\bea\label{GR}
D=e^{f_0}e^{-e_0}e^{f_0}e^{f_1}e^{-e_1}e^{f_1}.
\ena
It satisfies
\bea\label{DACT}
{\rm Ad}D\,E(z)=z^{-2}E(z),{\rm Ad}D\,H(z)=H(z)-2kz^{-1},
{\rm Ad}D\,F(z)=z^2F(z).
\ena
\begin{prop}
Consider the actions of $E(z)$ and $F(z)$ on the level $k$ representation
$\pi_j$.
We have
\bea
E(z)^k&=&c_j
e^{\int H_{<0}(z)}Dz^{H_0}e^{\int H_{>0}(z)},\label{EK}\\
F(z)^k&=&c_j
e^{-\int H_{<0}(z)}D^{-1}z^{-H_0}e^{-\int H_{>0}(z)}.\label{FK}
\ena
where $c_j=(-1)^jk!$.
\end{prop}
\proof
Note that the action of $D$ on $\pi_j$ is
determined (up to a constant multiple) by (\ref{DACT}).
For $k=1$ (\ref{EK}) and (\ref{FK}) follows from the well-known
result \cite{FK}. In the realization of \cite{FK}
the mutiplication by $e^{\alpha_1}$, where $\alpha_1$ is the simple root
of $sl_2$, is used instead of $D$.
It is easy to check (\ref{DACT}) for $e^{\alpha_1}$ in place of $D$.
Therefore, we have $e^{\alpha_1}=c_jD$ for some constant $c_j$.
By calculating $D|j\rangle\in\pi_j$ we obtain $c_j=(-1)^j$ for $k=1$.

Suppose we know that (\ref{EK}) is true for $k$.
Consider the tensor product of the representation of a level $k$
module with the currents $E(z)_1,H(z)_1,F(z)_1$, and a level $1$ module
with the currents $E(z)_2,H(z)_2,F(z)_2$. In the tensor product we
have a level $(k+1)$ action of $E(z)=E(z)_1+E(z)_2$. Because of (\ref{INT})
we have $E(z)^{k+1}=(k+1)E(z)_1^kE(z)_2$. Therefore $E(z)^{k+1}$ has
the same form as (\ref{EK}). The proof of (\ref{FK}) is similar.\qed

The operators $E(z)^k$ and $F(z)^k$ generate a vertex subalgebra.
Their operator product expansion reads as
\bea
E(z)^kF(w)^k=(z-w)^{-2k}\sum_{i\geq0}(z-w)^i:C^{(1)}_i(z):,
\ena
where each term $C^{(1)}_i(z)$ is a differential polynomial of $H(z)$
and, in particular, $C^{(1)}_0(z)$ is a constant.

Consider now the square root of the automorphism (\ref{DACT}):
\bea
U(E(z))=-z^{-1}E(z),U(F(z))=-zF(z),U(H(z))=H(z)-kz^{-1}.
\ena
This is an outer automorphism.

Let $\iota$ be the outer automorphhim of $\widehat{sl}_2$ induced from the
non-trivial Dynkin diagram automorphism. It is involutive and
$\iota(e_0)=e_1,\iota(h_0)=h_1,\iota(f_0)=f_1$. It acts on the currents:
\bea
\iota\bigl(E(z)\bigr)=zF(z),\iota\bigl(H(z)\bigr)=-H(z)+kz^{-1}.
\ena
It also acts from $\pi_j$ to $\pi_{k-j}$
in such a way that $\iota X\iota^{-1}=\iota(X)$ $(X\in\widehat{sl}_2)$
in ${\rm End}\left(\oplus_{j=0}^k\pi_j\right)$.

Set
\bea
D^{{1\over2}}=e^{f_0}e^{-e_0}e^{f_0}\iota.
\ena
Then we have
$\bigl(D^{{1\over2}}\bigr)^2=D$ and $D^{{1\over2}}XD^{-{1\over2}}=U(X)$
$(X\in\widehat{sl}_2)$ in ${\rm End}\left(\oplus_{j=0}^k\pi_j\right)$.

\begin{prop}
In ${\rm Hom}(\pi_j,\pi_{k-j})$ we have
\bea
\psi_0(z)&=&z^{{j\over2}}e^{{1\over2}\int H_{<0}(z)}
D^{{1\over2}}z^{{H_0\over2}}e^{{1\over2}\int H_{>0}(z)},\label{VERT}\\
\psi_k(z)&=&(-1)^{k-j}z^{{j\over2}}e^{-{1\over2}\int H_{<0}(z)}
D^{-{1\over2}}z^{-{H_0\over2}}e^{-{1\over2}\int H_{>0}(z)}.
\ena
\end{prop}
\proof
The operator $\psi_0(z):\pi_j\rightarrow\pi_{k-j}$
is determined (up to a constant multiple) by (\ref{DCOM}) and
the commutation relation
with $H(z)$ and $E(z)$. It is easy to check that the right hand side of
(\ref{VERT}) commutes with $E(z)$ and has the correct commutation relations
with $H(z)$. For $k=1$ one can also check directly (\ref{DCOM})
and (\ref{NOR}). The general case follows from the following coonsideration.

Suppose $k_l\in\Z_{>0}$ $(l=1,2)$ and consider
the representations $\left(\pi_{j_l}\right)_l$ of level $k_l$.
We put $(\phantom{W})_l$ only to distinguish the different values of level for
$l=1,2$. We use similar notations for the vertex operators
$\left(\psi_a(z)\right)_l$.

Consider the operator
\bea
\left(\psi_0(z)\right)_1\otimes\left(\psi_0(z)\right)_2:
\left(\pi_{j_1}\right)_1\otimes\left(\pi_{j_2}\right)_2\rightarrow
\left(\pi_{k_1-j_1}\right)_1\otimes\left(\pi_{k_2-j_2}\right)_2.
\ena
The algebra $\widehat{sl}_2$ of level $k_1+k_2$ is acting on
$\left(\pi_{j_1}\right)_1\otimes\left(\pi_{j_2}\right)_2$ and
$\left(\pi_{k_1-j_1}\right)_1\otimes\left(\pi_{k_2-j_2}\right)_2$
diagonally. Let us decompose 
\bea
\left(\pi_{j_1}\right)_1\otimes\left(\pi_{j_2}\right)_2=
\oplus_{\gamma=0}^{k_1+k_2}\pi_\gamma\otimes S_\gamma.
\ena
The representation
$\left(\pi_{k_1-j_1}\right)_1\otimes\left(\pi_{k_2-j_2}\right)_2$
is canonically isomorphic to
\bea
\oplus_{\gamma=0}^{k_1+k_2}\pi_{k_1+k_2-\gamma}\otimes S_\gamma.
\ena

From the explicit formula (\ref{VERT}) it follows that the
operator $\left(\psi_0(z)\right)_1\otimes\left(\psi_0(z)\right)_2$
is acting from
$\pi_\gamma\otimes S_\gamma$ to $\pi_{k_1+k_2-\gamma}\otimes S_\gamma$
and equal to $\psi_0(z)\otimes1$.\qed

Set $\phi_a(z)=z^{-j/2}\psi_a(z)$.
From (\ref{VERT}) we can deduce the operator product expansion:
\bea
\phi_0(z)\phi_0(w)=(z-w)^{{k\over2}}\left(\sum_{i\geq0}
:E(z)^kC^{(2)}_i(z):(z-w)^i\right),
\ena
where $C^{(2)}_i(z)$ is a differential polynomial of $H(z)$ and,
in particular, that $C^{(2)}_0$ is a constant.
A similar result holds if we replace $\phi_0,E$ with $\phi_k,F$. We have also
\bea
\phi_0(z)\phi_k(w)=(z-w)^{-{k\over2}}\left(\sum_{i\geq0}
:C^{(3)}_i(z):(z-w)^i\right),
\ena
where $C^{(3)}_i(z)$ is a differential polynomial of $H(z)$ and,
in particular, $C^{(3)}_0(z)$ is a constant.

Let $B$ be the Heisenberg algebra with the basis $\{a_j\}$ $(j\in\Z)$
and the relations $[a_i,a_j]=i\delta_{i+j,0}$. Let $\H_q$ be the irreducible
representation of $B$; $\H_q$ contains the vacuum vector $v(q)$ such that
$a_jv(q)=0$ $(j>0)$ and $a_0v(q)=qv(q)$. Introduce the vertex operator
$w(p,z)$:
\bea
w(p,z)=\exp\left(-p\sum_{i<0}{a_i\over i}z^{-i}\right)
T^pz^{pa_0}\exp\left(-p\sum_{i>0}{a_i\over i}z^{-i}\right).
\ena
The operator $T^p$ acts from $\H_q$ to $\H_{q+p}$.
It commutes with $a_i$ $(i\not=0)$ and satisfies $T^pv(q)=v(q+p)$.
The operator product expansion of $w(p,z)$ has the form
\bea
w(p_1,z_1)w(p_2,z_2)=(z_1-z_2)^{p_1p_2}\left(\sum_{i\geq0}
(z_1-z_2)^i S_i(z_1)\right)
\ena
where $S_0(z)=w(p_1+p_2,z)$. We have
\bea
[a_n,w(p,z)]=pz^nw(p,z).
\ena
Introduce now the operators acting from $\pi_j\otimes\H_q$
to $\pi_{k-j}\otimes\H_{q\pm\sqrt{k/2}}$
\bea
\varphi_a(z)&=&\phi_a(z)w\bigl(\sqrt{k/2},z\bigr),\\
\varphi^*_a(z)&=&(-1)^{k-j}\phi_a(z)w\bigl(-\sqrt{k/2},z\bigr).
\ena

\begin{dfn}
$A_k$ is the vertex operator algebra generated by
$\{\varphi_a(z)\},\{\varphi^*_a(z)\}$.
\end{dfn}

\subsection{Properties of $A_k$}
We now present some properties of the algebra $A_k$ without giving proofs.

Let $W\simeq V$ be the $(k+1)$-dimensional irreducible representation of
$gl_2$. The matrix $\pmatrix{1&0\cr0&1\cr}$
is acting on $W$ by the scalar $\sqrt{k\over2}$.
We write the basis of $W$ as $\varphi_a$ instead of $\psi_a$,
i.e., the linear map $\psi_a\mapsto\varphi_a$ is $sl_2$ linear.
Let $W^*$ be the dual space. We identify $W^*$ with $V$ by the invariant
coupling determined by $(\psi_0,\psi_k)=1$. We write the basis of $W^*$
as $\varphi^*_a$ instead of $\psi_a$. We denote the map $W\rightarrow W^*$,
$\varphi_a\mapsto\varphi^*_a$ by $*$.

The algebra $A_k$ is generated by the spaces
$W(z)$ with basis $\varphi_0(z),\ldots,\varphi_k(z)$ and
$W^*(z)$ with basis $\varphi^*_0(z),\ldots,\varphi^*_k(z)$.
For a vector $w\in W$ we will denote by $w(z)$ and $w^*(z)$ the corresponding
operators in $A_k$.

(a) Operators $w_1(z),w_2(z)$ $(w_i\in W)$ are commuting if
$k$ is even and skew-commuting if $k$ is odd. The same is true for
$w^*_1(z),w^*_2(z)$.

We call a vector $w\in W$ the highest if it is annihilated by
a nilpotent matrix $u$ in $gl_2$: $u(w)=0$. The set of all highest vectors form
a cone $K\subset W$. The set $K$ is an orbit by the $SL_2$ action on $W$.

(b) Let $w\in K$, then $w(z)w^{(l)}(z)=0$ for $l<k$ and
$w^*(z)w^{*(l)}(z)=0$ also for $l<k$.
Here $w^{(l)}(z)$ is the $l$-th derivative of $w(z)$.

(c) If $w_1,w_2\in W$, then
\bea
w_1(z_1)w^*_2(z_2)(z_1-z_2)^k=w^*_2(z_2)w_1(z_1)(z_2-z_1)^k.
\ena
The operator product expansion of $w_1(z_1)$ and $w^*_2(z_2)$ has the form:
\bea
w_1(z_1)w^*_2(z_2)=(z_1-z_2)^{-k}\left(S_0+(z_1-z_2)S_1(z_1)+\cdots\right)
\ena
where $S_0$ is a scalar, $S_0=(w_1,w_2)$, and
$S_1(z)$ is a linear combination of Heisenberg algebra
$a(z)=\sum a_iz^{-i}$ and $\widehat{sl}_2=\{E(z),H(z),F(z)\}$.
Altogether they constitute the algebra $\widehat{gl}_2$. Therefore, we see
that $A_k$ contains $\widehat{gl}_2$ as a Lie subalgebra. The algebra $A_k$ is generated by $\varphi_0(z),\varphi^*_0(z)$ and $\widehat{gl}_2$.

(d) Let $w\in K$. Then, the operator product of $w(z_1)$ and $w^*(z_2)$
has no singular terms. In particular, it means that
\bea
w(z_1)w^*(z_2)=(-1)^kw^*(z_2)w(z_1).
\ena
Cases $k=1$ and $k=2$. It is easy to see that for $k=1$ the operators
$\{\varphi_0(z)$, $\varphi_1(z)$, $\varphi^*_0(z)$, $\varphi^*_1(z)\}$
generate the usual Clifford algebra. In the case $k=2$ we have $6$ generators
$\{\varphi_0(z)$, $\varphi_1(z)$, $\varphi_2(z)$,
$\varphi^*_0(z)$, $\varphi^*_1(z)$, $\varphi^*_2(z)\}$.
The operator product of $\varphi_a(z_1)$ and $\varphi^*_b(z_2)$
starts from $(z_1-z_2)^{-2}$ with a constant coefficient. The next term is
$\widehat{gl}_2$. Therefore, all $\varphi_a(z),\varphi^*_a(z)$
generate (with respect to the bracket) the central extension of some Lie
algebra. It is $\widehat{sp}_4$ with level $1$.

(e) Fix two non-negative integers $k_1,k_2$. There is a homomorphism of
algebras,
\bea\label{MU}
\kappa:A_{k_1+k_2}\rightarrow A_{k_1}\otimes A_{k_2}.
\ena
The map $\kappa$ can be characterized by the following way. First we have
$\widehat{sl}_2$ in $A_{k_1+k_2}$. This $\widehat{sl}_2$ goes by a diagonal 
way in $\widehat{sl}_2\oplus\widehat{sl}_2$.
The same is true for the Heisenberg algebra $B$ in $A_{k_1+k_2}$.
Finally, $\varphi_0(z)\in A_{k_1+k_2}$ goes to the
$\bigl(\varphi_0(z)\bigr)_1\otimes\bigl(\varphi_0(z)\bigr)_2$,
and $\varphi^*_0(z)$ goes to $\bigl(\varphi^*_0(z)\bigr)_1
\otimes\bigl(\varphi^*_0(z)\bigr)_2$.

\subsection{Representations of $A_k$}

Let us form the following space:
\bea
\cdots\oplus\left(\pi_k\otimes \H_{-\sqrt{k\over2}}
\right)\oplus\left(\pi_0\otimes \H_0
\right)\oplus\left(\pi_k\otimes \H_{\sqrt{k\over2}}
\right)\oplus\left(\pi_0\otimes \H_{2\sqrt{k\over2}}\right)\oplus\cdots
\ena
It is clear that the operators $W(z)$ and $W^*(z)$
are acting on this space in such a way that they can be expanded in
$z^n$ $(n\in\Z)$. We can generalize this construction by writing
\bea
\left(\bigoplus_{s\in2\Z}\pi_j\otimes\H_{q+s\sqrt{k\over2}}\right)
\oplus\left(\bigoplus_{s\in2\Z+1}\pi_{k-j}\otimes
\H_{q+s\sqrt{k\over2}}\right)
\ena
where $q=m\sqrt{2\over k}$, $m=0,1,\ldots,k-1$ if $j$ is even,
and $q=(m+{1\over2})\sqrt{2\over k}$, $m=0,1,\ldots,k-1$ if $j$ is odd.

Therefore, our space is labeled by two numbers $j$ $(0\leq j\leq k)$
and $r=q\sqrt{2k}$, $(0\leq r\leq 2k-1)$ where $j+r$ is even.
Denote this space by $R(j,r)$. It is evident that $R(j,r)=R(j',r')$
if $j+j'=k$ and $r+k=r'\pmod{2k}$.

Without giving a proof we state
\begin{prop}
Each irreducible representation of the algebra $A_k$ has the form $R(j,r)$ for
some $j,r$. Therefore, the algebra $A_k$ has $k(k+1)/2$ irreducible
representations.
\end{prop}

Representations $R(j,r)$ of $A_k$ form the minimal models.
It means in particular, that representations $R(j,r)$ form the Verlinde
algebra. We describe it here. Denote by ${\cal V}_k$ the Verlinde algebra of $\widehat{sl}_2$ of level $k$. It is an algebra with
basis $\pi_0,\pi_1,\ldots,\pi_k$ and the product
\bea
\pi_i\pi_j=\pi_{i-j}+\pi_{i-j+2}+\cdots+\pi_s
\ena
where $j\leq i$ and $s=\min(i+j,2k-i-j)$.
Note that $\pi_0=1$.
Let ${\cal E}_{2k}$ be the algebra with basis
$1,\varepsilon,\varepsilon^2,\ldots,\varepsilon^{2k-1}$,
and product $\varepsilon^r\varepsilon^s
=\varepsilon^{r+s\pmod{2k}}$.
Consider the tensor product ${\cal V}_k\otimes{\cal E}_{2k}$. Define on
${\cal V}_k\otimes{\cal E}_{2k}$ the following operator $\gamma$:
\bea
\gamma(\pi_j\otimes\varepsilon^r)
=(-1)^{j+r}\pi_j\otimes\varepsilon^r.
\ena
Let $[{\cal V}_k\otimes{\cal E}_{2k}]^\gamma$ be the set of fixed elements
in ${\cal V}_k\otimes{\cal E}_{2k}$ by $\gamma$.
This is a subalgebra. The Verlinde algebra of $A_k$ is the quotient
$[{\cal V}_k\otimes{\cal E}_{2k}]^\gamma/J$ where $J$ is the ideal generated
by the relation $\pi_k\otimes\varepsilon^k=1$. The element
$\pi_j\otimes\varepsilon^r$
corresponds to the representation $R(j,r)$.

\subsection{Quadratic relations}
The operators
$W(z)$ acting on the irreducible representations of $A_k$ satisfy the
following quadratic relations.
\bea
&{\rm (a)}&\quad w_1(z)w_2(w)=(-1)^kw_2(w)w_1(z), w_1,w_2\in W,\\
&{\rm (b)}&\quad w(z)w^{(l)}(z)=0 (l<k)\hbox{ if $w\in K$}.\label{QR1}
\ena
We call them the relations (R).

Suppose that operators $W(z)$ are acting on some space $Q$ with the grading
$Q=Q_0\oplus Q_{-1}\oplus Q_{-2}\oplus\cdots$.
Choose vectors $\rho\in Q$ and $\rho^*\in Q^*$.
The correlation function is by definition the following matrix element:
\bea
\langle\rho^*,w_1(z_1)\cdots w_n(z_n)\rho\rangle, w_a\in W.
\ena
Changing $w_1,\ldots,w_n$ we get the vector-valued function
$F^*(\rho^*,\rho;z_1,\ldots,z_n)$ with values in
$\underbrace{W^*\otimes W^*\otimes\cdots \otimes W^*}_n$, which we will also
call the correlation function. It is evident that
$F^*(\rho^*,\rho;z_1,\ldots,z_n)$ is symmetric if $k$ is even,
and skew-symmetric if $k$ is odd. (We permute simultaneously
coordinates $z_i$ and components in $W^*\otimes W^*\otimes\cdots \otimes W^*$).
$F^*(\rho^*,\rho;z_1,\ldots,z_n)$ is a Laurent polynomial. $W^*$ is an
irreducible $sl_2$-module, so we can realize $W^*$ in the space of sections
of the line bundle on $\C P^1$. If we choose the coordinate $t$ on $\C P^1$
and trivialize the line bundle on $\C P^1$ then $W^*$ will be identified with
the space of polynomials in $t$ of degree less than or equal to $k$.
Therefore, the tensor product $W^*\otimes\cdots\otimes W^*$
can be identified with the space of polynomials in $t_1,\ldots,t_n$ of
degree less than or equal to $k$ in each variable $t_j$.
The correlation function $F^*(\rho^*,\rho;z_1,\ldots,z_n)$ can be viewed as
a scalar function in variable $X_1,\ldots,X_n$ where $X_j=(t_j,z_j)$;
we denote it by $F(\rho^*,\rho;X_1,\ldots,X_n)$.
\begin{prop}\label{P}
Suppose that operators $W(z)$satisfy the quadratic relations $(R)$.
Then the corresponding function $F(\rho^*,\rho,X_1,\ldots,X_n)$
has a zero of order at least $k$ if $X_i=X_j$ $1\leq i<j\leq n$.
In other words,
\bea
{\partial^{j_1}\over\partial t_i^{j_1}}
{\partial^{j_2}\over\partial z_i^{j_2}}
F(\rho^*,\rho,X_1,\ldots,X_n)\Big|_{X_i=X_j}=0
\hbox{ if $j_1+j_2<k-1$}.
\ena
\end{prop}
We prove this proposition in 2.3.

Consider the irreducible representation $R(j,r)$ of the algebra $A_k$.
We will call a vector in $\pi_i\otimes \H_q\subset R(j,r)$
``extremal'' if it is the product of an extremal vector in $\pi_i$ and the
highest weight vector $v(q)$.

\begin{prop}\label{AB}
Introduce the Fourier coefficients of $\varphi_a(z)$:
\bea
\varphi_a(z)=\sum_n\varphi_{a,n}z^{-n}.
\ena
Let $\omega\in R(j,r)$ be an extremal vector. Then
there exists a set of integers $N_a$ depending on $\omega$ and $a$
such that the following are valid.
\bea
&&(i)\quad\varphi_{a,n}\omega=0\hbox{  if $n>N_a$}.
\phantom{-------------------...}\label{HWC}
\ena

\quad$(ii)$\quad Let $B_\omega$ be the subspace of $R(j,r)$:
$B_\omega=\C[\varphi_{a,n};0\leq a\leq k,n\in\Z]\omega$
Then we have an isomorphism
\bea
\C[\varphi_{a,n};0\leq a\leq k,n\in\Z]/N
{\buildrel{\simeq}\over\rightarrow B_\omega}
\ena
where the ideal $N$ is generated by the relations $(R)$ and $(\ref{HWC})$.
\end{prop}

{\it Proof of $(i)$}
Using the automorphism of $A_k$ induced by $D^{1/2}$
we can reduce the proof to the case where $\omega=\omega_{j,j+2l}$:
\bea\label{W}
\omega_{j,j+2l}=|j\rangle\otimes v\Bigl({j+2l\over\sqrt{2k}}\Bigr).
\ena
We define the local energy function
$h:I_k\otimes I_k\rightarrow I_k$ where $I_k=\{0,1,\ldots,k\}$:
\bea
h_{a,b}=\min(a,k-b).\label{ENE}
\ena
For $\omega=\omega_{j,j+2l}$ we set
\bea
N_a=h_{a,k-j}-j-l.
\ena
Then we have
\bea
\varphi_{a,n}\omega_{j,j+2l}&=&0\hbox{ if $n>N_a$},\\
\varphi_{j,-l}\omega_{j,j+2l}&=&\omega_{k-j,j+2l+k}.\label{VAC}
\ena
These equalities follow by a direct calculation. Because of the automorphism
of $A_k$ induced by $T^{\sqrt{2/k}}$, 
it is enough to prove it for one value of $l$.\qed

The proof of $(ii)$ is given in 2.3.

\section{Monomial bases and correlation functions}
In the previous section we have introduced the vertex operator
algebra $A_k$ that are generated by $\varphi_a(z),\varphi^*_a(z)$ 
$(0\leq a\leq k)$.
In this section we construct a basis of its irreducible representations.
We interprete the quadratic relations of the vertex operators
as normal-ordering rules. The vectors in the basis are normal-ordered monomials
(see below for the precise definition).

\subsection{Symmetric and skew-symmetric tensors}
In this section we work on the Fourier coefficients of the vertex operators:
$\varphi_a(z)=\sum_{n\in\Z}\varphi_{a,n}z^{-n}$.
First we consider the $\varphi_{a,n}$ as abstract generators of an algebra $B_k$
where we assume only the (skew-) commutativity,
\bea
[\varphi_{a,m},\varphi_{b,n}]\mp=0\hbox{ if $k$ is even (odd).}
\ena
The $\widehat{sl}_2$-action is given by
\bea\label{ACTION}
e_0\varphi_a(z)=(a+1)z\varphi_{a+1}(z),
e_1\varphi_a(z)=(k-a+1)\varphi_{a-1}(z),\\
f_0\varphi_a(z)=(k-a+1)z^{-1}\varphi_{a-1}(z),
f_1\varphi_a(z)=(a+1)\varphi_{a+1}(z),\\
h_0\varphi_a(z)=-(n-2a)\varphi_a(z),
h_1\varphi_a(z)=(n-2a)\varphi_a(z).
\ena
The following diagram shows the action for $k=2$ schematically.
\newcommand{\mapdown}[1]{\Big\downarrow
  \rlap{$\vcenter{\hbox{$\scriptstyle#1\,$}}$}}
\newcommand{\mapnwarrow}[1]{\nwarrow\kern-5pt\raise5pt
  \hbox{$\vcenter{\hbox{$\scriptstyle#1\,$}}$}}

\[
\begin{array}{ccccccc}
 \cdots & \varphi_{0,-1} && \varphi_{0,0} && \varphi_{0,1} & \cdots \\
 & \mapdown{f_1} & \mapnwarrow{f_0} & \mapdown{f_1} & \mapnwarrow{f_0}
 & \mapdown{f_1} \\
 \cdots & \varphi_{1,-1} && \varphi_{1,0} && \varphi_{1,1} & \cdots \\
 & \mapdown{f_1} & \mapnwarrow{f_0} & \mapdown{f_1} & \mapnwarrow{f_0}
 & \mapdown{f_1} \\
 \cdots & \varphi_{2,-1} && \varphi_{2,0} && \varphi_{2,1} & \cdots 
\end{array}
\]

We introduce an ordering of the index set $\{0,1,\ldots,k\}\times\Z$:
$(a,m)<(b,n)$ if and only if $m<n$, or $m=n$ and $a>b$.
Note that $f_0$ and $f_1$ are lowering operators in this ordering.

If $n$ is even, we have $B_k=\oplus_{s=0}^\infty B^{(s)}_k$ where
\bea
B^{(s)}_k=\oplus_{(a_1,n_1)\leq\cdots\leq(a_s,n_s)}
\C\varphi_{a_1,n_1}\cdots\varphi_{a_s,n_s}.
\ena
If $n$ is odd, we have
\bea
B^{(s)}_k=\oplus_{(a_1,n_1)<\cdots<(a_s,n_s)}
\C\varphi_{a_1,n_1}\cdots\varphi_{a_s,n_s}.
\ena

In order to handle the quadratic relations (\ref{QR1}), we need a completion of
the algebra. For $N\in\Z$, let $B_{k,N}$ be the ideal of $B_k$ that is
generated by the set of elements $\{\varphi_{a,m};m>N\}$. We set
$B^{(s)}_{k,N}=B^{(s)}_k\cap B_{k,N}$, and define the
completion of the algebra $B_k$.
\bea
\bar B_k=\oplus_{s=0}^\infty\bar B^{(s)}_k,\\
\bar B^{(s)}_k=\lim_{\buildrel\leftarrow\over n}B^{(s)}_k/B^{(s)}_{k,n}.
\ena
The $\widehat{sl}_2$-action extends to $\bar B_k$.

Let $U\subset W\otimes W$ be a subspace. We define $D^{(m)}U$ to be the
subspace of $\bar B^{(2)}_k$ spanned by the Fourier coefficients of
\bea\label{DMW}
D^{(m)}w(z)=\sum_{a,b}c_{a,b}
\varphi_a(z)\varphi^{(m)}_b(z)
\ena
where
$w=\sum_{a,b}c_{a,b}\varphi_a\otimes\varphi_b\in U$
and $\varphi_b^{(m)}(z)$ is the $m$-th derivative of $\varphi_b(z)$.

We have the decomposition
\bea
W\otimes W=S\oplus A
\ena
where $S$ and $A$ are the symmetric and the skew-symmetric tensors,
respectively. They are invariant with respect to the $sl_2$ action.
In fact, we have the irreducible decomposition of $S$ and $A$ as follows.
\bea
S&=&\C^{2k+1}\oplus\C^{2k-3}\oplus\cdots\nonumber\\
&=&\cases{\C^1\oplus\C^5\oplus\cdots\hbox{ if $k$ is even};\cr
\C^3\oplus\C^7\oplus\cdots\hbox{ if $k$ is odd},\cr}\\
A&=&\C^{2k-1}\oplus\C^{2k-5}\oplus\cdots\nonumber\\
&=&\cases{\C^3\oplus\C^7\oplus\cdots\hbox{ if $k$ is even};\cr
\C^1\oplus\C^5\oplus\cdots\hbox{ if $k$ is odd}.\cr}
\ena

Set $\Omega=\sum_{m=0}^\infty D^{(m)}(W\otimes W)$.
The space of the quadratic relations is an $\widehat{sl}_2$-invariant
subspace of $\Omega$. We now determine these spaces (see Propositions
\ref{P1} and \ref{R}).

In abuse of notation we write $D^{(m)}\C^{2j+1}$ $(0\leq j\leq k)$
for $D^{(m)}U$ where
$U$ is the unique $(2j+1)$-dimensional component of $W\otimes W$.
We have
\bea
D^{(0)}(W\otimes W)=
D^{(0)}\C^1\oplus D^{(0)}\C^5\oplus\cdots.
\ena
If $m\geq1$ we have
\bea
D^{(m)}(W\otimes W)=
D^{(m)}\C^1\oplus D^{(m)}\C^3\oplus\cdots.
\ena
These are direct sums. However, the sum of
$D^{(m)}(W\otimes W)$ for $m=0,1,\ldots$ is not direct.
In fact, we have
\begin{prop}\label{P1}
$(i)$ We have the irreducible  decomposition of the $sl_2$ module $\Omega$:
\bea\label{DECOM}
\Omega&=&(D^{(0)}\C^1\oplus D^{(0)}\C^5\oplus\cdots)\nonumber\\
&&\oplus(D^{(1)}\C^3\oplus D^{(1)}\C^7\oplus\cdots)\nonumber\\
&&\oplus(D^{(2)}\C^1\oplus D^{(2)}\C^5\oplus\cdots)\nonumber\\
&&\oplus(D^{(3)}\C^3\oplus D^{(3)}\C^7\oplus\cdots)\oplus\cdots.
\ena
$(ii)$ The space $D^{(m)}\C^{2j+1}$ which does not appear in the above
decomposition is contained in $\sum_{l=0}^{m-1}D^{(l)}\C^{2j+1}$.
\end{prop}
\proof
The statement $(ii)$ follows from the (skew-)commutativity
of $\varphi_{a,n}$ and the following equality.
\bea
\varphi^{(m)}_a(z)\otimes\varphi_b(z)
-(-1)^m\varphi_a(z)\otimes\varphi^{(m)}_b(z)\nonumber\\
=\sum_{j=0}^{m-1}(-1)^j\left(m\atop j\right)
{d^{m-j}\over dz^{m-j}}\left(\varphi_a(z)
\otimes{d^j\over dz^j}\varphi_b(z)\right).
\ena
We now prove that the summands in (\ref{DECOM}) are linearly independent.
Suppose $w=\sum c_{ab}\varphi_a\otimes\varphi_b\in S$ 
and $k+m$ is even, or $w\in A$ and $k+m$ is odd.
Consider the degree $d$ term of $D^{(m)}w$:
\bea
&&\sum_{d_1+d_2=d}c_{ab}
(-d_2)(-d_2-1)\cdots(-d_2-m+1)\varphi_{a,d_1}\varphi_{b,d_2}\nonumber\\
&=&\sum_{d_1+d_2=d}{1\over2}c_{ab}
\{(-d_2)(-d_2-1)\cdots(-d_2-m+1)\nonumber\\
&&+(-1)^m(-d_1)(-d_1-1)\cdots(-d_1-m+1)\}
\varphi_{a,d_1}\varphi_{b,d_2}\nonumber\\
&=&\sum_{2d_1\leq d}c_{ab}
(\epsilon(d_1)d_1^m+\cdots)\varphi_{a,d_1}\varphi_{b,d-d_1},
\ena
where
\bea
\epsilon(d_1)=\cases{1\hbox{ if $2d_1=d$};\cr 2\hbox{ otherwise}.}
\ena
Namely, the coefficients increase in the $m$-th power of $d_1$.
The proof of $(i)$ is over.\qed

The $sl_2$ structure of $\Omega$ is clear from (\ref{DECOM}). We now determine
the $\widehat{sl}_2$ structure of $\Omega$. This is not semi-simple.
The following table illustrates the closure relations of the
$\widehat{sl}_2$-action.
\bea\label{TABLE}
\matrix{&m=0&&m=1&&m=2&&m=3&&m=4\cr
k=0&1&&&&1&&&&1\cr
&&\nwarrow&&\swarrow&&\nwarrow&&\swarrow&\cr
k=1&&&3&&&&3&&\cr
&&\swarrow&&\nwarrow&&\swarrow&&\nwarrow&\cr
k=2&{\bf5}&&&&5&&&&5\cr
&&\nwarrow&&\swarrow&&\nwarrow&&\swarrow&\cr
k=3&&&{\bf7}&&&&7&&\cr
&&\swarrow&&\nwarrow&&\swarrow&&\nwarrow&\cr
k=4&{\bf9}&&&&{\bf9}&&&&9\cr
&&\nwarrow&&\swarrow&&\nwarrow&&\swarrow&\cr
k=5&&&{\bf11}&&&&{\bf11}&&\cr
&&\swarrow&&\nwarrow&&\swarrow&&\nwarrow&\cr
k=6&{\bf13}&&&&{\bf13}&&&&{\bf13}\cr}
\ena
The integer $2j+1$ in the $m$-th column signifies $D^{(m)}\C^{2j+1}$
in (\ref{DECOM}). We have
\begin{prop}\label{P2}
$(i)$ Let $\Omega_1\subset\Omega$ be an $\widehat{sl}_2$-invariant subspace.
If $D^{(m)}\C^{2j+1}\subset\Omega_1$, then
$D^{(m-1)}\C^{2j-1}\oplus D^{(m-1)}\C^{2j+3}\subset\Omega_1$.

$(ii)$ Consider a subdiagram $I$ of $(\ref{TABLE})$.
Let $\Omega(I)$ be the union of the corresponding subspaces in $\Omega$.
$\Omega(I)$ is $\widehat{sl}_2$-invariant if and only if $I$ is closed
with respect to the arrows.
\end{prop}

\proof
The `only if' part of $(ii)$ follows from $(i)$. Let us show $(i)$ and the
`if' part of $(ii)$. Suppose that $D^{(m)}\C^{2j+1}\subset\Omega(I)$. Let
$v^{(j)}_0=\sum c^{(j)}_{ab}\varphi_a\otimes\varphi_b
\in\C^{2j+1}\subset W\otimes W$ be the highest weight vector.
We will show that the Fourier coefficients of $e_0D^{(m)}v^{(j)}_0(z)$
(see (\ref{DMW}))
belongs to the closure of $D^{(m)}\C^{2j+1}$ in the sense of the arrows.
Because $[e_0,f_1]=0$, the same statement for 
$e_0\left(D^{(m)}f_1^lv_0^{(j)}(z)\right)$ then follows.

We have
\bea\label{E0}
&&e_0D^{(m)}v^{(j)}_0(z)
=zf_1\left(D^{(m)}v^{(j)}_0(z)\right)+v_1(z),\nonumber\\
&&v_1(z)=m\sum c^{(j)}_{ab}
\varphi_a(z)(f_1\varphi_b(z))^{(m-1)}.
\ena

It is easy to see that $(1\otimes f_1)v^{(j)}_0$ is a linear combination of
$f_1^2v^{(j+1)}_0$, $f_1v^{(j)}_0$ and $v^{(j-1)}_0$ with non-zero
coefficients. The proposition follows from this.\qed

\subsection{Quadratic relations and normal-ordering rules}
We assume that $k\geq2$. We denote by $Q^{(2)}_k$ the
$\widehat{sl}_2$-invariant subspace of $\bar B^{(2)}_k$ that is generated by
$\varphi_0(z)\varphi_0^{(k-2)}(z)$ (or equivalently by
$\varphi_k(z)\varphi_k^{(k-2)}(z)$).

For example, if $k=2$, the $\widehat{sl}_2$ action generates
\bea
\varphi_0(z)^2,\varphi_0(z)\varphi_1(z),2\varphi_0(z)\varphi_1(z)
+\varphi_1(z)^2,\varphi_1(z)\varphi_2(z),\varphi_2(z)^2.
\ena
The following proposition follows immediately from
Propositions \ref{P1} and \ref{P2}.
\begin{prop}\label{R}
The $sl_2$ decomposition of $Q^{(2)}_k$ is given by
\bea\label{RDEC}
Q^{(2)}_k=\oplus_{j\equiv m(\bmod 2),\atop2m+5\leq2j+1\leq2k+1}
D^{(m)}\C^{2j+1}
\ena
\end{prop}

The summands in (\ref{RDEC}) are indicated in (\ref{TABLE}) by the boldface
letters. For example, $Q^{(2)}_2$ consists of $D^{(0)}\C^5$.
$Q^{(2)}_3$ of $D^{(0)}\C^5$ and $D^{(1)}\C^7$,
$Q^{(2)}_4$ of $D^{(0)}\C^5$,
$D^{(1)}\C^7$, $D^{(0)}\C^9$, $D^{(2)}\C^9$, etc.
We define the algebra $R_k$ to be the quotient of $\bar B_k$ by the ideal
generated by $Q^{(2)}_k$.

Our next aim is to rewrite the quadratic relations as the normal-ordering
rules. Consider first a special case, $k=2$ and the quadratic relation
$\varphi_0(z)^2=0$. In terms of the Fourier coefficients the relations read as
\bea
\sum_{j_1+j_2=j}\varphi_{0,j_1}\varphi_{0,j_2}=0.
\ena
We can rewrite them as follows:
\bea
\varphi_{0,l}\varphi_{0,l}=
-2\sum_{n=1}^\infty\varphi_{0,l-n}\varphi_{0,l+n},\\
\varphi_{0,l}\varphi_{0,l+1}=
-\sum_{n=1}^\infty\varphi_{0,l-n}\varphi_{0,l+1+n}.
\ena
The product $\varphi_{0,j_1}\varphi_{0,j_2}$
where $(j_1,j_2)$ satisfying $j_1\geq j_2-1$
are written as linear combinations of products
$\varphi_{0,l_1}\varphi_{0,l_2}$ where $(l_1,l_2)$
satisfying $l_1\leq l_2-2$. We call such relations normal-ordering rules.

We now consider the general situation.
The way to rewrite the quadratic relations as normal-ordering rules is
not unique. Set $J=\{0,1,\ldots,k\}\times\Z$.
Choose and fix a subset $O\subset J\otimes J$.
A pair $((a,m),(b,n))$ is called ``normal-ordered''
if $((a,m),(b,n))\in O$. A product $\varphi_{a_1,\j_1}\cdots\varphi_{a_s,j_s}$
is called normal-ordered if each pair
$((a_l,j_l),(a_{l+1},j_{l+1}))$
$(1\leq l\leq s-1)$ is normal-ordered.

A set of normal-ordered pairs $O$ is called good if the following are valid.

(a) Any product $\varphi_{a_1,j_1}\cdots\varphi_{a_s,j_s}$
can be rewritten as a linear combination of normal-ordered products
by using the normal-ordering rules.

(b) Normal-ordered products are linearly independent.\\
Our aim is to find a good subset $O$.

Define the local energy function $h:J\times J\rightarrow\Z$ by
\bea
h((a,m),(b,n))=h_{a,b}-m+n
\ena
where $h_{a,b}$ is given by (\ref{ENE}).
We set
\bea
O=\{((a,m),(b,n));h((a,m),(b,n))\geq k\}.
\ena
Our goal is to prove
\begin{prop}
The set of normal-ordered pairs $O$ given above is good.
\end{prop}
We give a proof of (a) in this section, and (b) in 2.3.
For the proof of (a) it is enough to show the following
\begin{prop}If $\varphi_{a_0,m_0}\varphi_{b_0,n_0}$ is not normal-ordered,
one can rewrite it as a linear combination of
$\varphi_{a,m}\varphi_{b,n}$
where $(a,m)<(a_0,m_0)$.
\end{prop}

\proof
If $(a_0,m_0)>(b_0,n_0)$, or if $k$ is odd and
$(a_0,m_0)=(b_0,n_0)$,
then the assertion is obvious because of the (skew-)commutativity
of $\varphi_{a,n}$. Therefore, we assume that
\bea\label{NA}
(a_0,m_0)\leq(b_0,n_0)\hbox{ if $k$ is even,}\quad
(a_0,m_0)<(b_0,n_0)\hbox{ if $k$ is odd.}
\ena
For $m_0\leq n_0$ we set
\bea
F_{m_0,n_0}=\{(a_0,b_0);\Bigl((a_0,m_0),(b_0,n_0)\Bigr)\not
\in O\hbox{ and (\ref{NA}) is satisfied}\},
\ena
and
$\left(F_{m_0,n_0}\right)_d=\{(a,b)\in F_{m_0,n_0};a+b=d\}$.

Let $U^*_{m_0,n_0}\left(\hbox{resp. }(U^*_{m_0,n_0})_d\right)
\subset\left(W\otimes W\right)^*$
be the subspace spanned by the set of elements
$\varphi^*_{a_0}\otimes\varphi^*_{b_0}$ where
$(a_0,b_0)\in F_{m_0,n_0}\left(\hbox{resp. }(F_{m_0,n_0})_d\right)$.
Here
$\langle\varphi^*_{a_1}\otimes\varphi^*_{b_1},
\varphi_{a_2}\otimes\varphi_{b_2}\rangle
=\delta_{a_1a_2}\delta_{b_1b_2}$.
Note that there is an action of $sl_2$ on 
$\left(W\otimes W\right)^*$ such that
$\langle Xw^*,w\rangle+\langle w^*,Xw\rangle=0$ for
$w^*\in\left(W\otimes W\right)^*$
and
$w\in W\otimes W$.

For $w=\sum c_{ab}\varphi_a\otimes\varphi_b
\in W\otimes W$ we set
$w_{m,n}=\sum c_{ab}\varphi_{a,m}\varphi_{b,n}\in R_k$.
We denote the map $w\rightarrow w_{m,n}$ by $\pi_{m,n}$.

Suppose that $w\in\C^{2j+1}\subset W\otimes W$.
If $n_0-m_0\geq0$ and $D^{(n_0-m_0)}\C^{2j+1}$ belongs to (\ref{RDEC}),
then by using Proposition \ref{P1} $(ii)$ we have
\bea
\sum c_{ab}\varphi_a(z)\left(z^{n_0-1}
\varphi_b(z)\right)^{(n_0-m_0)}=0
\ena
in $R_k$. The degree $m_0+n_0$ part of this equation contains
$(-1)^{n_0-m_0}(n_0-m_0)!w_{m_0,n_0}$, and the rest is a linear combination
of $w_{m_0-l,n_0+l}$ $(l\geq1)$. Similarly,
if $n_0-m_0\geq1$ and $D^{(n_0-m_0-1)}\C^{2j+1}$ belongs
to (\ref{RDEC}), we have 
\bea
\sum c_{ab}\varphi_a(z)\left(z^{n_0-1}
\varphi_b(z)\right)^{(n_0-m_0-1)}=0,
\ena
in $R_k$. The degree $m_0+n_0$ part of this equation contains
$(-1)^{n_0-m_0-1}2(n_0-m_0-1)!w_{m_0,n_0}$,
and the rest is a linear combination of $w_{m_0-l,n_0+l}$ $(l\geq1)$.
Therefore, we have

$(i)$ if $m_0=n_0$ and $w\in\C^{4j+1}\subset 
W\otimes W$ $(1\leq j\leq[{k\over2}])$ then $w_{m_0,n_0}$
is equal to a linear combination of $w_{m_0-l,n_0+l}$ $(l\geq1)$;

$(ii)$ if $n_0-m_0\geq1$ and $ w\in\C^{2j+1}\subset 
W\otimes W$ $(n_0-m_0+1\leq j\leq k)$ then $w_{m_0,n_0}$
is equal to a linear combination of $w_{m_0-l,n_0+l}$ $(l\geq1)$.

Let $U_{m_0,n_0}\subset
W\otimes W$ be the sum of the irreducible components appearing
in $(i)$ or $(ii)$.

Now we prove the statement of the proposition for
$(a_0,b_0)\in F_{m_0,n_0}$. We assume that
$0\leq a_0+b_0\leq k$. The case $k\leq a_0+b_0\leq2k$
is similar.

For $U\subset W\otimes W$ set
$U_d=\{w\in U;h_1w=2(k-d)w\}$.
We claim that the the coupling between $(U^*_{m_0,n_0})_d$ and
$(U_{m_0,n_0})_d$ is non-degenerate. If $m_0=n_0$ and $0\leq d\leq k-1$,
this is clear because
\bea
(U_{m_0,n_0})_d=\cases{S_d&if $k$ is even;\cr
A_d&if $k$ is odd.\cr}
\ena
It is also clear if $m_0<n_0$ and
$0\leq d\leq k-n_0+m_0-1$ because $(U_{m_0,n_0})_d=(W\otimes W)_d$.
We have, in particular, that if $0\leq a_0+b_0\leq k-n_0+m_0-1$
then $\varphi_{a_0,m_0}\varphi_{b_0,n_0}\in\pi_{m_0n_0}(U_{m_0,n_0})$,
and therefore, $\varphi_{a_0,m_0}\varphi_{b_0,n_0}$
can be written as a linear combination of
$\varphi_{a_0,m_0-l}\varphi_{b_0,n_0+l}$ $(l\geq1)$.

The non-degeneracy of the coupling for $k-n_0+m_0\leq d\leq k$
reduces to the case $d=k-n_0+m_0-1$ because
$f_1^{d-k+n_0-m_0+1}:(U_{m_0,n_0})_{k-n_0+m_0-1}\rightarrow(U_{m_0,n_0})_d$
and
$f_1^{d-k+n_0-m_0+1}:(U^*_{m_0,n_0})_d\rightarrow(U^*_{m_0,n_0})_{k-n_0+m_0-1}$
are both isomorphisms of vector spaces.

Thus, we have shown that if $(a_0,b_0)\in F_{m_0,n_0}$
the product $\varphi_{a_0,m_0}\varphi_{b_0,n_0}$
can be written as a linear combination of
$\varphi_{a,m}\varphi_{b,n}$ where, in addition to
$a+b=a_0+b_0, m+n=m_0+n_0$, we have
``$m<m_0,n>n_0$'' or ``$m=m_0,n=n_0$ and $(a,b)\not\in F_{m_0,n_0}$''.

Observe that in the latter case we have
$(a,m_0)<(a_0,m_0)$. It will finish the proof.\qed

\subsection{Correlation functions}
We use the realization of the dual space,
$W^*=\oplus_{a=0}^k\C t^a$
where $\langle t^a,\varphi_b\rangle=\delta_{a,b}$.
The dual $sl_2$ action is given by
\bea\label{SL2}
e_1=t^2{d\over dt}-kt,h_1=2t{d\over dt}-k,f_1=-{d\over dt}.
\ena
In this picture $(W^*)^{\otimes s}$ is nothing but the space of
polynomials in $t_1,\ldots,t_s$ whose degree in each $t_j$ is less
than or equal to $k$.

Fix $0\leq r\leq k$.
We extend this realization to the dual space of $B^{(s)}_k/N^{(s)}_{k,r}$,
where $N^{(s)}_{k,r}$ is the subspace of $B^{(s)}_k$ that is the union
of $B^{(s-1)}_k\varphi_{a,n}$ for $n>\min(a-r,0)$.
The dual space
$\left(B^{(s)}_k/N^{(s)}_{k,r}\right)^*$ is realized as
the space of polynomials
in $t_1,\ldots,t_s,z_1,\ldots,z_s$ satisfying the following conditions:

$(i)$ $f(t_1,\ldots,t_s,z_1,\ldots,z_s)$ is symmetric if $k$ is even, and
skew-symmetric if $k$ is odd, with respect to the permutation of
$(t_1,z_1),\ldots,(t_s,z_s)$;

$(ii)$ the degree of $f(t_1,\ldots,t_s,z_1,\ldots,z_s)$ in each $t_j$ is less
than or equal to $k$.

$(iii)$ 
\bea
{\partial^{j_1}\over\partial t_1^{j_1}}
{\partial^{j_2}\over\partial z_1^{j_2}}
f(t_1,\ldots,t_s,z_1,\ldots,z_s)|_{t_1=z_1=0}=0
\hbox{ for $j_1+j_2<r$.}
\ena

The dual coupling is induced from
\bea
\langle \prod_{j=1}^st_j^{a_j}\prod_{j=1}^sz_j^{m_j},
\varphi_{b_1,n_1}\otimes\cdots\otimes\varphi_{b_s,n_s}\rangle
=\delta_{a_1,b_1}\cdots\delta_{a_s,b_s}
\delta_{m_1,-n_1}\cdots\delta_{m_s,-n_s}.
\ena

Now we consider the quotient ${\cal H}^{(s)}_{k,r}$ of
$B^{(s)}_k/N^{(s)}_{k,r}$
by the quadratic relations discussed in 2.2.
The dual space ${\cal H}^{*(s)}_{k,r}$ is a subspace of
$\left(B^{(s)}_k/N^{(s)}_{k,r}\right)^*$.

A polynomial $f(t_1,\ldots,t_s,z_1,\ldots,z_s)\in
\left(B^{(s)}_k/N^{(s)}_{k,r}\right)^*$ 
belongs to ${\cal H}^{*(s)}_{k,r}$ if and only if
the following condition is satisfied:

$(iv)$ if $\C^{2j+1}\subset W\otimes W$ is such that
$D^{(m)}\C^{2j+1}$ is a component of (\ref{RDEC})
then we have
\bea\label{QR}
\langle {\partial^m\over\partial z_1^m}f(t_1,\ldots,t_s,z_1,\ldots,z_s),
\C^{2j+1}\rangle|_{z_1=z_2}=0.
\ena

We call the function $f(t_1,\ldots,t_s,z_1,\ldots,z_s)$ satisfying
$(i),(ii),(iii),(iv)$ an $s$-particle correlation function.

\begin{prop}
The condition $(\ref{QR})$ is equivalent to
\bea\label{37}
{\partial^{j_1}\over\partial t_1^{j_1}}
{\partial^{j_2}\over\partial z_1^{j_2}}
f(t_1,\ldots,t_s,z_1,\ldots,z_s)|_{t_1=t_2,z_1=z_2}=0
\hbox{ for $j_1+j_2<k$.}
\ena
\end{prop}
\proof
We write
\bea
{\partial^m\over\partial z_1^m}f(t_1,\ldots,t_s,z_1,\ldots,z_s)|_{z_1=z_2}
=\sum_i(u_i)_{12}\otimes(v_i)_{3\cdots s}
\subset W^*\otimes\cdots\otimes W^*
\ena
where $u_i\in W^*\otimes W^*$. Then, the condition (\ref{QR})
is equivalent to
$$u_i\in\oplus_{j=0}^m{\bf C}^{2j+1}\subset W^*\otimes W^*.$$
Using (\ref{SL2}) one can show that this is further equivalent to the condition
that $u_i$ has the factor $(t_1-t_2)^{k-m}$.
The statement of the proposition follows from this.\qed

The proof of Proposition \ref{P} is similar to this.

\begin{prop}\label{PR}
The set of the normal-ordered monomials in ${\cal H}^{(s)}_{k,r}$, i.e.,
the set of vectors of the form
\bea
\varphi_{a_1,n_1}\cdots\varphi_{a_s,n_s}\quad
\ena
where $h((a_j,n_j),(a_{j+1},n_{j+1}))\geq k$
for all $1\leq j\leq s-1$ and $n_s\leq\min(a_s-r,0)$,
is linearly independent.
\end{prop}
\proof
Set
\bea
I^{(s)}_{k,r}&=&\{((a_1,n_1),\ldots,(a_s,n_s));
h((a_j,n_j),(a_{j+1},n_{j+1}))\geq k\nonumber\\
&&\hbox{ for all $1\leq j\leq s-1$ and $n_s\leq\min(a_s-r,0)$}\}.\nonumber\\
\ena
We define an order in $I^{(s)}_{k,r}$ by saying that
\bea
&&((a_1,m_1),\ldots,(a_s,m_s))
<((b_1,n_1),\ldots,(b_s,n_s))\nonumber\\
&&\hbox{ if and only if for some $l$ we have }\nonumber\\
&&(a_j,m_j)=(b_j,n_j)\hbox{ for $1\leq j\leq l-1$}
\hbox{ and }(a_l,m_l)<(b_l,n_l).
\ena

For $\rho=((b_1,n_1)\ldots,(b_s,n_s))\in I^{(s)}_{k,r}$ we set
$v_\rho=\varphi_{b_1,n_1}\cdots\varphi_{b_s,n_s}
\in{\cal H}^{(s)}_{k,r}$.
We will construct a set of correlation functions
$\{f_\kappa(t_1,\ldots,t_s,z_1,\ldots,z_s);\kappa\in I^{(s)}_{k,r}\}$
such that
\bea\label{TRI}
\langle
f_\kappa(t_1,\ldots,t_s,z_1,\ldots,z_s),v_\rho\rangle=
\cases{1& if $\kappa=\rho$;\cr
0& if $\kappa>\rho$.\cr}
\ena

For $0\leq a\leq k$ and $m\leq0$ we define a $k$-component
monomial in $t$ and $z$:
\bea\label{43}
(P^{(1)}_{a,m},\ldots,P^{(k)}_{a,m})
=(t^{a^{(1)}}z^{m^{(1)}},\ldots,t^{a^{(k)}}z^{m^{(k)}})
\ena
where $a^{(1)},\ldots,a^{(k)}$ and $m^{(1)},\ldots,m^{(k)}$
are uniquely determined by the following conditions:
\bea
&&\sum_ja^{(j)}=a,\quad\sum_jm^{(j)}=-m,\\
&&m^{(1)}\geq m^{(2)}\geq\cdots\geq m^{(k)}\geq m^{(1)}-1,\\
&&m^{(1)}+a^{(1)}\geq m^{(2)}+a^{(2)}\geq\cdots
m^{(k)}+a^{(k)}\geq m^{(1)}+a^{(1)}-1.
\ena

For $p_1,\ldots,p_s\in\C[t,z]$ we define
\bea
p_1\odot\cdots\odot p_s
=\sum_{\sigma\in S_s}{\rm sgn}\,\sigma\ p_1(t_{\sigma_1},z_{\sigma_1})\cdots 
p_s(t_{\sigma_s},z_{\sigma_s}).
\ena
For $\kappa=((a_1,m_1)\ldots,(a_s,m_s))\in I^{(s)}_{k,r}$ we define
\bea
f_\kappa(t_1,\ldots,t_s,z_1,\ldots,z_s)
=\prod_{j=1}^kP^{(j)}_{a_1,m_1}\odot\cdots\odot P^{(j)}_{a_s,m_s}.
\ena
One can easily check $(i),(ii),(iii)$,(\ref{37}) and (\ref{TRI}). \qed

Finally we give a proof of Proposition \ref{AB} $(ii)$. The proof is similar
to the above proof of Proposition \ref{PR}. Since the surjectivity is clear,
it is enough to show the injectivity. Without loss of generality we assume
that $w=\omega_{r,r}$ $(0\leq r\leq k)$.

It is enough to show that the set of vectors
\bea
\{v_\rho\omega_{r,r};\rho\in I^{(s)}_{k,r}\}\label{50}
\ena
is linearly independent.

We reduce the proof to the level $1$ case by using the algebra map
$A_k\rightarrow\underbrace{A_1\otimes\cdots\otimes A_1}_k$
(see (\ref{MU})). Let us write $\omega^{(k)}_{a,r}$ and
$\varphi^{(k)}_a(z)$ for $\omega_{a,r}$ and $\varphi_a(z)$
of level $k$.
We realize $\omega^{(k)}_{r,r}$
as ${\underbrace{\omega^{(1)}_{1,1}\otimes\cdots\otimes\omega^{(1)}_{1,1}}_r}
\otimes{\underbrace{\omega^{(1)}_{0,0}\otimes\cdots\otimes
\omega^{(1)}_{0,0}}_{k-r}}$. The action of $\varphi^{(k)}_a(z)$
is realized as
\bea
a!\sum_{a_1+\cdots+a_k=a}
\varphi^{(1)}_{a_1}(z)\otimes\cdots\otimes\varphi^{(1)}_{a_k}(z).
\ena
We will show the linear independence of the vectors (\ref{50})
in this realization.

Consider the vectors
\bea
\omega^i_{a_1,\ldots,a_s;m_1,\ldots,m_s}=
\varphi^{(1)}_{a_1,m_1}\cdots\varphi^{(1)}_{a_s,m_s}
\omega^{(1)}_{i,i}
\ena
such that $-2m_1+a_1>\cdots>-2m_s+a_s\geq i$ for $i=0,1$, and their dual 
vectors $\omega^{i*}_{a_1,\ldots,a_s;m_1,\ldots,m_s}$.

For $0\leq a\leq k$, $m\leq0$ and $1\leq j\leq k$, we define
$A^{(j)}_{a,m}\in\{0,1\}$ and $M^{(j)}_{a,m}\in\Z_{\leq0}$ by
\bea
A^{(j)}_{a,m}=a^{(j)},\quad M^{(j)}_{a,m}=-m^{(j)},
\ena
where $a^{(j)},m^{(j)}$ are given in (\ref{43}).

For $\kappa=((a_1,m_1),\ldots,(a_s,m_s))
\in I^{(s)}_{k,r}$ and $1\leq j\leq k$ we set
\bea
\omega^*_j(\kappa)=\omega^{i*}_{A^{(j)}_{a_1,m_1},\ldots,A^{(j)}_{a_s,m_s};
M^{(j)}_{a_1,m_1},\ldots,M^{(j)}_{a_s,m_s}},
\ena
where $i=\cases{1& if $1\leq j\leq r$;\cr0& if $r+1\leq j\leq k$.\cr}$
Then, it is easy to see that
\bea
\langle\omega^*_1(\kappa)\otimes\cdots\otimes
\omega^*_k(\kappa),v_\rho\omega^{(k)}_{r,r}\rangle=
\cases{c_\rho& if $\kappa=\rho$;\cr0& if $\kappa>\rho$,\cr}
\ena
where $c_\rho$ is a non-zero constant.\qed

From (\ref{VAC}) the vector $\omega_{r,r+2l}\in\pi_r
\otimes{\cal H}_{{r+2l\over\sqrt{2k}}}$ can be formally written as
\bea
\omega_{r,r+2l}
=\varphi_{k-r,k-r-l}\varphi_{r,k-l}\varphi_{k-r,2k-r-l}\varphi_{r,2k-l}\cdots.
\ena
For $j\geq1$ we define
\bea
a^{(r,r+2l)}_j&=&\cases{k-r& if $j$ is odd;\cr r& if $j$ is even,\cr}\\
n^{(r,r+2l)}_j&=&\cases{{(j+1)k\over2}-r-l& if $j$ is odd;\cr
{jk\over2}-l& if $j$ is even.\cr}\ena
We call a sequence $p=(a_j,n_j)_{j\geq1}$ a {\it path} which belongs to
$\omega_{r,r+2l}$ if the following are satisfied.
\bea
&(i)&0\leq a_j\leq k\quad n_j\in\Z,\\
&(ii)&a_j=a^{(r,r+2l)}_j,\quad n_j=n^{(r,r+2l)}_j\hbox{ if $j>\hskip-3pt>0$},\\
&(iii)&((a_j,n_j),(a_{j+1},n_{j+1}))\hbox{ is normal-ordered}.
\ena
We denote the set of the paths which belong to $\omega_{r,r+2l}$
by ${\cal P}_{r,r+2l}$. We can associate a vector
$\omega_p\in\pi_r\otimes{\cal H}_{{r+2l\over\sqrt{2k}}}$
with each path $p\in{\cal P}_{r,r+2l}$:
\bea
\omega_p=\varphi_{a_1,n_1}\varphi_{a_2,n_2}\cdots.
\ena
From what we have proved it follows that the vectors $\omega_p$ are
linearly independent. The weight of $\omega_p$ is given by the formula.
\bea
{\rm wt}(\omega_p)-{\rm wt}(\omega_{r,r+2l})
=-\sum_{j=1}^\infty(a_j-a^{(r,r+2l)}_j)\alpha_1
+\sum_{j=1}^\infty(n_j-n^{(r,r+2l)}_j)\delta.
\ena

We call a path
$p=(a_j,n_j)_{j\geq1}\in{\cal P}_{r,r+2l}$
``reduced'' if
\bea
h((a_j,n_j),(a_{j+1},n_{j+1}))=k
\ena
for all $j\geq1$. We denote the set of the reduced paths in
${\cal P}_{r,r+2l}$ by ${\cal P}^{\rm red}_{r,r+2l}$.

We conclude the paper by the following

\begin{prop}
The space $\pi_r\otimes{\cal H}_{{r+2l\over\sqrt{2k}}}$
is spanned by the set of vectors
$\omega_p=\varphi_{a_1,n_1}\varphi_{a_2,n_2}\cdots$
for $p\in{\cal P}_{r,r+2l}$.
\end{prop}
\proof
We compare the characters of the space
$\pi_r\otimes{\cal H}_{{r+2l\over\sqrt{2k}}}$,
and its subspace spanned by the above vectors.
Let us denote the former by $\chi_{r,r+2l}$
and the latter $\chi({\cal P}_{r,r+2l})$.
We also denote the character of the space spanned by the vectors corresponding
to the reduced paths
by $\chi({\cal P}^{{\rm red}}_{r,r+2l})$. Then we have
\bea
\chi({\cal P}_{r,r+2l})=
{\chi({\cal P}^{{\rm red}}_{r,r+2l})\over\prod_{n=1}^\infty(1-e^{-n\delta})}.
\ena
For a reduced path $p=(a_j,n_j)_{j\geq1}$ we have
\bea
\sum_{j=1}^\infty(n_j-n^{(r,r+2l)}_j)
=\sum_{j=1}^\infty j(h_{a_j,a_{j+1}}-h_{a^{(r,r+2l)}_j,a^{(r,r+2l)}_j}).
\ena
Therefore, the assertion $\chi_{r,r+2l}=\chi({\cal P}_{r,r+2l})$
follows from the known fact, Theorem 1.2 of \cite{DJKMO}.\qed

\end{document}